\RequirePackage[leqno]{amsmath}
\documentclass[leqno]{amsart}
\usepackage{amsmath}
\makeatletter
\usepackage[utf8]{inputenc}     
\usepackage{multicol}
\usepackage{color, graphics}
\usepackage[demo]{graphicx}
\usepackage{efbox,graphicx}
\efboxsetup{linecolor=gray,linewidth=0.8pt}
\usepackage{floatrow}
\usepackage{pgfplots}
\usepackage{graphicx, multicol, latexsym, amsmath, amssymb}
\usepackage{pgfplots}
\usepackage{tikz}
\usepackage{amsfonts}
\usepackage{enumitem}
\setlist[itemize]{topsep=0pt,after=\vspace{1.5\baselineskip}}
 \usepackage[colorlinks=true]{hyperref}
\usepackage{empheq}
\usepackage[T1]{fontenc}
\usepackage{tabularx}
\usepackage[leftcaption]{sidecap}
\sidecaptionvpos{figure}{m}
\usepackage{tikz}
\usepackage{subfigure}
\usepackage[margin=0.5in]{geometry}
\usetikzlibrary{backgrounds,automata}

\newcount\rc@count
\let\rc@clearconstantlist\empty
\newcommand\rc@clearconstant[1]{\global\expandafter\let\csname rc@const@#1\endcsname\undefined}
\newcommand\resetconstants[1]{%
    \def\rc@constname{#1}
    \global\rc@count=1\relax 
    \bgroup 
        \let\\\rc@clearconstant 
        \rc@clearconstantlist
        \global\let\rc@clearconstantlist\empty 
    \egroup
}
\newcommand\const[1]{%
    \@ifundefined{rc@const@#1}{%
        \expandafter\xdef\csname rc@const@#1\endcsname{%
           \noexpand\rc@useconst{\rc@constname}{\the\rc@count}%
        }%
        \g@addto@macro\rc@clearconstantlist{\\{\mathrm{#1}}}%
        \global\advance\rc@count1\relax
    }{}%
    \csname rc@const@#1\endcsname
}
\newcommand\rc@useconst[2]{{#1}\textsubscript{#2}}
\setlist[itemize]{noitemsep, topsep=0pt}

\def\R{\mathbb R} \def\N{\mathbb N}

\def\R{\mathbb R} \def\N{\mathbb N} 
\def\TM{T_{max}} 
\def
\@cite
#1#2{[#1\if@tempswa, #2\fi]}
\makeatother
\newtheorem{theorem}{Theorem}[section]
\newtheorem{corollary}[theorem]{Corollary}
\newtheorem{example}[theorem]{Example}
\newtheorem{lemma}[theorem]{Lemma}

\newtheorem{definition}[theorem]{Definition}
\newtheorem{remark}{Remark}

\newcounter{cnstcnt}

\title[Boundedness in a Keller--Segel system with production and consumption] 
{Combining effects ensuring boundedness in an attraction-repulsion chemotaxis model with production and consumption}
\author[Tongxing Li, Silvia Frassu and Giuseppe Viglialoro]{}
\subjclass[2020]{Primary: 35A01, 35K55, 35Q92. Secondary:  92C17.}
\keywords{Chemotaxis, Global existence, Boundedness, Nonlinear production, Consumption, Logistic source. \\
\textit{$^*$Corresponding author}: giuseppe.viglialoro@unica.it}

\begin{document}

\maketitle

\centerline{\scshape{\scshape{Tongxing Li$^{\natural}$, Silvia Frassu$^{\sharp}$ \and Giuseppe Viglialoro$^{\sharp,*}$}}}
\medskip
{
\medskip
\centerline{$^{\natural}$School of Control Science and Engineering} 
\centerline{Shandong University}
\centerline{Jinan, Shandong, 250061 (P. R. China)}
\medskip
}
{
\medskip
\centerline{$^\sharp$Dipartimento di Matematica e Informatica}
\centerline{Universit\`{a} di Cagliari}
\centerline{Via Ospedale 72, 09124. Cagliari (Italy)}
\medskip
}
\bigskip
\begin{abstract}
This paper is framed in a series of studies on attraction-repulsion chemotaxis models combining different effects: nonlinear diffusion and sensitivities and logistic sources, for the dynamics of the cell density, and consumption and/or production impacts, for those of the chemicals. In particular, herein we focus on the situation where the signal responsible of gathering tendencies for the particles' distribution is produced, while the opposite counterpart is consumed. In such a sense, this research complements the results in \cite{FrassuLiViglialoro} and \cite{ChiyoFrassuViglialoro-Att-Rep-2022}, where the chemicals evolve according to different laws.
\end{abstract}
\resetconstants{c}
\section{Introduction, motivations and state of the art}\label{Intro}
\subsection{The attractive and the repulsive models}\label{IntroDiscussionSection}
During the last 50 years, precisely since the advent of basic biological models idealizing chemotaxis phenomena, many related variants have attracted the interest of the mathematical community. Indeed, several technical difficulties encouraging and stimulating researchers in the field arise in such problems.

For our purposes, we are interested in a model resulting in a combination of this aggregative signal-production mechanism  
\begin{equation}\label{problemOriginalKS} 
u_t= \Delta u - \chi \nabla \cdot (u \nabla v) =\nabla \cdot(\nabla u-\chi u\nabla v)\quad \textrm{and} \quad 
\tau v_t=\Delta v-v+u, \quad \textrm{ in } \Omega \times (0,\TM),
\end{equation}
and this repulsive signal-consumption one 
\begin{equation}\label{problemOriginalKSCosnumption}
u_t= \Delta u+ \xi \nabla \cdot (u \nabla w)=\nabla \cdot(\nabla u+\xi u\nabla w) \quad \textrm{and} \quad 
 w_t=\Delta w-u w, \quad \textrm{ in } \Omega \times (0,\TM).
\end{equation}
Herein $\Omega$ is a bounded and smooth domain of $\R^n$, with
$n\in \N$, $\tau \in \{0,1\}$, $\chi, \xi>0$ and $\TM>0$. Moreover, to properly interpret these formulations in the context of real biological phenomena, we also specify the following: $u=u(x,t)$ denotes a particular cell distribution (populations, organisms) at position $x\in \Omega$ and time $t\in (0,\TM)$, and $v=v(x,t)$ and $w=w(x,t)$ stand for the concentrations of chemoattractant and chemorepellent (i.e., a chemical signal that causes cells to attract and repel each other, respectively). In this way, system
\eqref{problemOriginalKS}, equipped with homogeneous Neumann boundary conditions, describes the dynamics of a cells' density in an isolated region whose (attractive) signal is produced by the cells themselves; in particular, the motion is influenced by the flux $1  \nabla u-\chi u \nabla v$, which idealizes how a constant (precisely, $1$) diffusive action is contrasted by an attractive one ($-\chi u$). Oppositely,  in model \eqref{problemOriginalKSCosnumption} the (repulsive) signal is absorbed by the cells, which are perturbed by the flux $1  \nabla u+\xi u \nabla w$ involving the same diffusion rate but a repulsion action ($+\xi u$). Naturally, apart from the boundary conditions, for the well-posedness of these systems some initial configurations for the cell density and the chemical signals have to be established; we will identify $u(x,0),\tau v(x,0)$ and $ w(x,0)$ respectively with $u_0=u_0(x), \tau v_0=\tau v_0(x)$ and $ w_0= w_0(x)$. 
\begin{remark}
According to the pioneering Keller--Segel models (\cite{K-S-1970,Keller-1971-MC,Keller-1971-TBC}), 
in problem \eqref{problemOriginalKS} the cell density $u$ evolves under the attractive action of the chemical $v$ (for instance a nutrient), which is secreted by the same cells. In parallel, if a chemical is absorbed by the cells (for instance when bacteria move toward the oxygen concentration), essentially $-v+u$ is replaced by $-vu$. Nevertheless, and even though the analysis of this paper is merely theoretical, it is worthwhile emphasizing that the signal-consumption model given in \eqref{problemOriginalKSCosnumption} for repulsive substances (for instance toxins) has as well some biological interpretations; indeed, \textit{phagocytes} are cells that protect the body by ingesting harmful foreign particles, bacteria, and dead or dying cells. For instance, during the phagocytosis process, the hepatic cells filter toxic substances, engulf them  and convert them into harmless substances or make sure they are released into the surrounding environment. (See \cite[$\S$1 and $\S$2]{MauriceHallett-MolecularBook}.)
\end{remark}
Coming back to our mathematical problems, as far as \eqref{problemOriginalKS} is concerned, it is known that if the attractive signal $v$ increases with $u$, the natural spread of the cells' density could suffer a gathering process and very high and spatially concentrated spike formations could appear. This mechanism, called \textit{chemotactic collapse} (or \textit{blow-up at finite time}), is connected to $\chi$, $m=\int_\Omega u_0(x)dx$ (the initial mass of the particle distribution) and $n$. Mathematically, $\TM$ is finite and the solution $(u,v)$ is local and unbounded at $\TM$. In this regard, if for $n=1$ blow-up phenomena are excluded (and in this case $\TM=\infty$ and $(u,v)$ is global),  when $n\geq 2$ whenever $m \chi$  surpasses a certain critical value $m_\chi$, the system might present the aforementioned chemotactic collapse, whereas for lower values of $m\chi$ respect to $m_\chi$ no instability appears in the motion of the cells. These results are included in more general analyses, dealing with the existence and properties (globality, uniform boundedness, or blow-up) of solutions to the Cauchy problem associated to \eqref{problemOriginalKS}, specially for the parabolic-elliptic version corresponding to the choice $\tau=0$; for details we mention \cite{OsYagUnidim,HerreroVelazquez,JaLu,Nagai,WinklAggre}, and references therein. On the other hand, for nonlinear production models, when in problem \eqref{problemOriginalKS} the linear segregation $f(u)=u$ behaves as $u^k$, with  $0<k<\frac{2}{n}$ ($n\geq1$), uniform boundedness of all its solutions is established in \cite{LiuTaoFullyParNonlinearProd}. Conversely, for a further simplification of the parabolic-elliptic formulation,  in spatially radial contexts it is known  (see \cite{WinklerNoNLinearanalysisSublinearProduction}) that  the value $\frac{2}{n}$ is critical, in the sense that the boundedness of solutions is guaranteed for any $n\geq 1$ and $0<k<\frac{2}{n}$, whereas for $k>\frac{2}{n}$ blow-up phenomena may occur. But, in the case of linear production rate, there is even more: The combination of \eqref{problemOriginalKS} with terms describing population growth or decay, for example logistic sources (see \cite{verhulst}), is very natural. In this situation, the equation for the particles' density  reads $u_t = \Delta u -\chi \nabla \cdot (u \nabla v) + f(u)$, where $f(u)$ may have the form $f(u)=\lambda u- \mu u^{\beta}$, $\lambda,\mu>0$ and $\beta>1$, and mathematical intuition indicates that the presence of the superlinear dampening effect should benefit the boundedness of solutions. Actually, this has only been shown for large values of $\mu$  (if $\beta=2$, see \cite{TelloWinkParEl}, \cite{W0}), whereas the possibility of blow-up was established for certain values of $\beta>1$, first for dimensions $5$ or higher \cite{WinDespiteLogistic}, when the equation for $v$ reads $0=\Delta v-\frac{1}{|\Omega|}\int_\Omega u(x,t)dx+u$ (see also \cite{FuestCriticalNoDEA} for a recent improvement of \cite{WinDespiteLogistic}), but later in \cite{Winkler_ZAMP-FiniteTimeLowDimension} also in 3-dimensional domains (with $0=\Delta v-v +u$ and $\beta<\frac{7}{6}$).

For the case where the repulsive signal $w$ is produced, which corresponds to replace in \eqref{problemOriginalKSCosnumption} the product $-uw$ with the sum 
$-u+w$, to our knowledge no result on the blow-up scenario is available and this is conceivable due to the repulsive nature of the mechanism; moreover, the corresponding state of the art is somehow poor  (see, for instance, \cite{Mock74SIAM,Mock75JMAA} for analyses on similar contexts). Oppositely, as we now will discuss, rather rich is the level of knowledge for attraction-repulsion
chemotaxis problems with linear or nonlinear production.
\subsection{The attraction-repulsion model: combining effects}\label{SectionAttr-Repulsion}
The general form of an attraction-repulsion chemotaxis model with logistic source and linear or nonlinear production can be formulated as 
\begin{equation}\label{problemAttRep}
	 \hspace*{-0.2cm}
	\begin{cases}
		u_t= \nabla \cdot \left(D(u)\nabla u - S(u)\nabla v + T(u)
		\nabla w\right)+h(u)  & \text{ in } \Omega \times (0,\TM),\\
		\tau_1 v_t= \Delta v +\phi(u,v)  & \text{ in } \Omega \times (0,\TM),\\
		\tau_2 w_t= \Delta w +\psi(u,w)& \text{ in } \Omega \times (0,\TM),\\
		u_{\nu}=v_{\nu}=w_{\nu}=0 & \text{ on } \partial \Omega \times (0,\TM),\\
		u(x,0)=u_0(x),\; \tau_1 v(x,0)=\tau_1 v_0(x),\; \tau_2 w(x,0)=\tau_2 w_0(x)& x \in \bar\Omega,\\ 
	\end{cases}
\end{equation}
where $D=D(u), S=S(u), T=T(u), h(u)$ are functions with a certain regularity and proper behavior, and $\phi=\phi(u,v)\simeq -v+u^a$ and
$\psi=\psi(u,w)\simeq-w+ u^b$ for some $a,b>0$, and $\tau_1,\tau_2\in\{0,1\}$. Moreover, further regular initial data  $u_0(x)\geq 0$, $ \tau_1 v_0(x)\geq 0$ and $ \tau_2 w_0(x)\geq 0$ are as well given, the subscript $\nu$ in $(\cdot)_\nu$ indicates the outward normal derivative on $\partial \Omega$, whereas $\TM \in (0,\infty]$ the maximal instant of time up to which solutions to the system do exist. 

These models have practical interest due to their real applicability to the inflammation observed in Alzheimer's disease where microglia secrete favorable attractant as well as unfavorable repellent; in particular, in the case of linear diffusion, sensitivities and productions (i.e., $D=S=T\equiv 1$ and $a=b=1$), and in absence of dampening external terms (i.e., $h\equiv 0$), \cite{Luca2003Alzheimer} deals with the description of the gathering mechanisms for \eqref{problemAttRep} and with dimensional, numerical and experimental analyses, in bounded intervals for the case $\tau_1=\tau_2=0$. 

Further, from the stricter mathematical point of view, when $D(u)\simeq u^{m_1}, S(u)\simeq u^{m_2}, T(u)\simeq u^{m_3}$ and $h(u)\simeq \lambda u- \mu u^\beta$ (where $m_1,m_2,m_3,\lambda,\mu,\beta$ attain some real values), we mention that for linear and nonlinear productions, criteria toward boundedness, long time behaviors and  blow-up issues for related solutions to \eqref{problemAttRep}  when $\tau_1,\tau_2\in \{0,1\}$ can be found in 
\cite{GuoJiangZhengAttr-Rep,LI-LiAttrRepuls,TaoWanM3ASAttrRep,VIGLIALORO-JMAA-BlowUp-Attr-Rep,YUGUOZHENG-Attr-Repul,ViglialoroMatNacAttr-Repul,LiangEtAlAtt-RepNonLinProdLogist-2020,XinluEtAl2022-Asymp-AttRepNonlinProd,ChiyoMarrasTanakaYokota2021,GuoqiangBin-2022-3DAttRep,GuoqiangBinATT-RepNonlinDiffSensLogistic}.
\section{The model and the main claims}\label{IntroSection}
\subsection{The model within a general frame}
Having in mind model \eqref{problemAttRep}, our research fits in a general project where the equations for the chemical signals $v$ and $w$ are both of consumption type, or one of production type and the other absorptive, or viceversa. Indeed, as indicated in the above set of papers, to our knowledge the actual state of the art concerns situations where the chemoattractant and the chemorepellent are both produced. 

So, let us give some details concerning the  results obtained so far within the aforementioned project. Specifically, since in this article we will focus on the situation where the signal responsible of gathering tendencies for the particles' distribution is produced, while the opposite counterpart is consumed, this research complements these already available results: \cite{FrassuLiViglialoro} and \cite{ChiyoFrassuViglialoro-Att-Rep-2022}. In the former the chemoattractant is consumed and the chemorepellent segregated, whilst in the latter both are absorbed. Specifically, for $D(u)\simeq u^{m_1}, S(u)\simeq u^{m_2}, T(u)\simeq u^{m_3}$ and $h(u)\simeq \lambda u- \mu u^\beta$ boundedness properties of solutions to \eqref{problemAttRep} are obtained  
\begin{enumerate}[label=\Roman*)]
\item \label{ItemImprovment} (\cite{FrassuLiViglialoro}) for $\tau_1=1$, $\tau_2=0$, $\phi(u,v)\simeq  -u^\alpha v$, $\psi(u,w)\simeq -\delta w+u^l$, $\delta>0$, $\alpha\in (0,1]$, $l\geq 1$, 
\item  (\cite{ChiyoFrassuViglialoro-Att-Rep-2022}) for $\tau_1=\tau_2=1$, $\phi(u,v)\simeq -u^\alpha v$, $\psi(u,w)\simeq -u^\gamma w$, $\alpha,\gamma \in (0,1]$
\end{enumerate}
whenever, according to the range of $\alpha$ and $l$ and  $\alpha$ and $\gamma$ respectively, $m_1$ is sufficiently larger than a function depending on the other remaining parameters of the problem. In particular, biological intuitions would suggest that the model where $v$ tends to vanish and $w$ to increase appears more inclined to provide boundedness with respect to the one with opposite effects of the chemicals, so that a natural question we want to deal with is:
\begin{enumerate}[label=$\mathcal{Q}_\arabic*$]
\item\label{Question}: If the assumptions in item \ref{ItemImprovment} are replaced with $\tau_1=0$, $\tau_2=1$, $\phi(u,v)\simeq -\delta v+ u^l$, $\psi(u,w)\simeq -u^\alpha w$, $\alpha\in (0,1]$, $l\geq 1$, are effectively the related conditions toward boundedness less sharp than those obtained in \cite{FrassuLiViglialoro}?
\end{enumerate}
We will analyze this question in Remark \ref{RemarkComparison} below. 
\subsection{Presentation of the Theorems} 
Conforming to all of the above, we study the following attraction-repulsion chemotaxis model
\begin{equation}\label{problem}
\begin{cases}
u_t= \nabla \cdot ((u+1)^{m_1-1}\nabla u - \chi u(u+1)^{m_2-1}\nabla v+\xi u(u+1)^{m_3-1}\nabla w) + h(u) & \text{ in } \Omega \times (0,\TM),\\
0= \Delta v - \delta v + g(u)& \text{ in } \Omega \times (0,\TM),\\
w_t=\Delta w-f(u)w  & \text{ in } \Omega \times (0,\TM),\\
u_{\nu}=v_{\nu}=w_{\nu}=0 & \text{ on } \partial \Omega \times (0,\TM),\\
u(x,0)=u_0(x), \; w(x,0)=w_0(x) & x \in \bar\Omega,
\end{cases}
\end{equation}
where $\Omega \subset \R^n$, $n \geq 2$, is a bounded and smooth domain, $\chi,\xi,\delta>0$, $m_1,m_2,m_3\in\R$, $f(u), g(u)$ and $h(u)$ are reasonably regular functions generalizing the prototypes $f(u)=K u^\alpha$, $g(u)=\gamma u^l$, and $h(u)=k u - \mu u^{\beta}$ with $K,\gamma, \mu>0$, $k \in \R$ and suitable $\alpha, l, \beta>0$. Once nonnegative initial configurations $u_0$ and $w_0$ are fixed, we aim at deriving sufficient conditions involving the above data so to ensure that problem \eqref{problem}
admits classical solutions which are global and uniformly bounded in time. 

To this purpose, we first give this 
\begin{definition}\label{ClassicalSolutionDefi}
We say that $(u,v,w)=(u(x,t),v(x,t),w(x,t))$ defined for $(x,t) \in \bar{\Omega}\times [0,\TM)$ is a classical solution to problem \eqref{problem}, if 
\begin{equation*}
0\leq u,w\in C^0(\bar{\Omega}\times [0,\TM))\cap  C^{2,1}(\bar{\Omega}\times (0,\TM)), 0\leq v\in C^0(\bar{\Omega}\times [0,\TM))\cap  C^{2,0}(\bar{\Omega}\times (0,\TM)),
\end{equation*}
and $u,v,w$ pointwisely satisfy all the relations in \eqref{problem}.
\end{definition}
Moreover, we require that $f$, $g$ and $h$ fulfill 
\begin{equation}\label{f}
f,g \in C^1(\R) \quad \textrm{with} \quad   0\leq f(s)\leq Ks^{\alpha}  \textrm{ and } \gamma s^l\leq g(s)\leq \gamma s(s+1)^{l-1},\;  \textrm{for some}\; K,\,\gamma,\,\alpha>0,\, l\geq 1 \quad \textrm{and all } s \geq 0,
\end{equation}
and 
\begin{equation}\label{h}
h \in C^1(\R) \quad \textrm{with} \quad h(0)\geq 0  \textrm{ and } h(s)\leq k s-\mu s^{\beta}, \quad  \textrm{for some}\quad k \in \R,\,\mu>0,\, \beta>1\, \quad \textrm{and all } s \geq 0.
\end{equation}
Then we establish these two theorems.
\begin{theorem}[The non-logistic model]\label{MainTheorem}
Let $\Omega$ be a smooth and bounded domain of $\mathbb{R}^n$, with $n\geq 2$, $\chi, \xi, \delta$ positive reals. Moreover, for  $m_1, m_2, m_3 \in \R$ and $h\equiv 0$, let $f$ and $g$ fulfill \eqref{f} for each of the following cases:
\begin{enumerate}[label=$C_{\roman*}$)]
\item \label{A3} $\alpha \in \left(0, \frac{1}{n}\right]$ and $m_1>m_2+l-\frac{2}{n}$ and $m_1>\min\left\{2m_3-(m_2+l),\max\left\{2m_3-1,\frac{n-2}{n}\right\}, m_3 - \frac{1}{n}\right\}$,
\item \label{A4} $\alpha \in \left(\frac{1}{n},\frac{2}{n}\right)$ and $m_1>m_2+l-\frac{2}{n}$ and $m_1>\min\left\{2m_3+1-(m_2+l),\max\left\{2m_3,\frac{n-2}{n}\right\}, m_3 + \alpha - \frac{2}{n}\right\}$,
\item \label{A5} $\alpha \in \left[\frac{2}{n},1\right]$ and $m_1>m_2+l-\frac{2}{n}$ and $m_1>\min\left\{2m_3+1-(m_2+l),\max\left\{2m_3,\frac{n-2}{n}\right\}, m_3 + \frac{n\alpha-2}{n\alpha-1}\right\}$.
\end{enumerate}
Then for any initial data $(u_0,w_0)\in (W^{1,\infty}(\Omega))^2$, with $u_0, w_0\geq 0$ on $\bar{\Omega}$,  problem \eqref{problem} admits a unique solution (in the sense of Definition \ref{ClassicalSolutionDefi}), for which $\TM=\infty$ and $(u,v,w) \in (L^\infty((0, \infty);L^{\infty}(\Omega)))^3.$

The same statement holds if we replace \ref{A5} with the condition 
\begin{enumerate}[label=$C'_{\roman*}$)]
\setcounter{enumi}{2}
\item \label{Ap5} $ \alpha \in \left[\frac{2}{n},1\right] \text{ and } m_2>2-l \text{ and } m_1>m_2+l-\frac{2}{n} \text{ and } m_1>\min\left\{2m_3+1-(m_2+l),\max\left\{2m_3,\frac{n-2}{n}\right\}\right\}.$
\end{enumerate}
\end{theorem}
\begin{theorem}[The logistic model]\label{MainTheorem1}
Under the same assumptions of Theorem \ref{MainTheorem}, let $f$, $g$ and $h$ fulfill \eqref{f} and \eqref{h} for each of the following cases:
\begin{enumerate}[label=$C_{\roman*}$)]
\setcounter{enumi}{3}
\item \label{A1} $\alpha \in \left(0, \frac{1}{n}\right]$ and $m_2<\beta-l$ and $m_1>\min\left\{2m_3-\beta, \max\left\{2m_3-1,\frac{n-2}{n}\right\}\right\}$,
\item \label{A2} $\alpha \in \left(\frac{1}{n}, 1\right)$ and $m_2<\beta-l$ and $m_1>\min\left\{2m_3-\beta+1, \max\left\{2m_3,\frac{n-2}{n}\right\}\right\}$.
\end{enumerate}
Then the same claim holds true.  In addition, if $\beta>2$ the conclusion applies even for $\alpha=1$.
\end{theorem}
By adaptations of the above results, we can also prove this 
\begin{corollary}[The linear diffusion and sensitivities case for the logistic model: $m_1=m_2=m_3=1$]\label{MainCorollary1}
Under the same assumptions of Theorem \ref{MainTheorem1}, let $m_1=m_2=m_3=1$. Moreover, for $\alpha \in (0,1]$ and $\beta>1+l$, let $f$, $g$ and $h$ satisfy \eqref{f} and \eqref{h}, respectively; then the same claim holds true. 
Conversely, for $\beta=1+l$ the following assumptions
\begin{enumerate}[label=$\tilde{C}_{\roman*}$)]
\item \label{A1B} $\alpha \in \left(0, \frac{1}{n}\right]$, $l\geq 1$ and $\mu>2^{l} \chi \gamma \frac{n-2}{n}=:\mathcal{L}(n)$,
\item \label{A2B} $\alpha \in \left(\frac{1}{n}, 1\right]$ and $l>1$ and $\mu>\mathcal{L}(n)$,
\item \label{A3B} $\alpha \in \left(\frac{1}{n}, 1\right)$ and $l=1$ and $\mu>2 \chi \gamma \frac{n-2}{n}+ \xi^2 2^{\frac{4-n}{n}} \left(\frac{n-2}{n+2}\right)^{\frac{n+2}{n}} n^{\frac{n+2}{n}} (n+1)^{\frac{2}{n}} \|\xi w_0\|_{L^{\infty}(\Omega)}^{\frac{4}{n}}=:\mathcal{L}(n)+\mathcal{M}(n)$,
\item \label{A4B} $\alpha=l=1$ and $\mu>\mathcal{L}(n)+\mathcal{M}(n) + 2^{\frac{n+2}{2}} K^{\frac{n+2}{2}} \left(\frac{3n-2}{n+2}\right)^{\frac{n+2}{4}} ((n-2)(n^2+n))^{\frac{n-2}{4}}  \|w_0\|_{L^{\infty}(\Omega)}^n=:\mathcal{L}(n)+\mathcal{M}(n)+\mathcal{N}(n)$,
\end{enumerate}
must be accomplished to obtain the claim.
\end{corollary}
\begin{remark}\label{RemarkOnLinearCase}
Let us give some considerations:
\begin{itemize}
\item [$\diamond$] From Corollary \ref{MainCorollary1} we deduce that the last condition is the strongest one; in particular it is obtained for the limit values of $\alpha$ and $l$. As soon as such values leave the limits, the conditions weaken, particularly for small values of $\alpha$. We also observe that $\mathcal{L}(2)=\mathcal{M}(2)=0$ and $\mathcal{N}(2)=4 K^2\|w_0\|^2_{L^\infty(\Omega)}$, so that in 2-dimensional settings the same corollary does not require restriction on $\mu$ except for the limit case $\alpha=l=1$.
\item [$\diamond$]  In general, Theorem \ref{MainTheorem} does not hold in the linear diffusion and sensitivities case; indeed, by substituting $m_1=m_2=m_3=1$ relation $m_1>m_2+l-\frac{2}{n}$ would require $l<\frac{2}{n}$, contradicting our assumption $l\geq1$. Conversely, the machinery works (see Remark \ref{JustificationRemark} at the end of the paper) for $\alpha \in (0,1]$, $n=2$ and $l=1$, but under the constrain $\chi \gamma<\frac{4}{C_{GN}}$, being $C_{GN}$ a constant depending on the domain. (This scenario sensitively changes with respect Corollary \ref{MainCorollary1}, where the presence of the logistic source allows us to develop a more exhaustive analysis also in higher dimensions.)
\end{itemize}
\end{remark}
\begin{remark}[On the question \ref{Question}: effects of produced/consumed chemoattractant/chemorepellent and viceversa]\label{RemarkComparison} This investigation and \cite{FrassuLiViglialoro} have in common the equation for the cells' density but opposite situations for the dynamics of the chemoattractant and the chemorepellent; in the present article the first is produced with rate $u^l$, $l\geq 1$, and the second  absorbed proportionally to $u^\alpha$, $\alpha \in (0,1]$, in the other the first is consumed with rate $u^\alpha$, $\alpha \in (0,1]$ and the second secreted proportionally to $u^l$, $l\geq 1.$ In view of this parallelism, question \ref{Question} above appears meaningful; nevertheless, as we can see in what follows, from the strict mathematical point of view the answer is not positive in any scenario.

First, in order to compare the results derived in Theorem \ref{MainTheorem} with those obtained in \cite{FrassuLiViglialoro}, by taking into account \cite[Remark 3]{FrassuLiViglialoro}, let us observe that we can rephrase assumptions $A_{i}$), $A_{ii})$ and $A_{iii})$ of \cite[Theorem 2.1]{FrassuLiViglialoro} in a more complete way, and more precisely as: 
\begin{itemize}
\item [] $\mathcal{A}_i$) $\alpha \in \left(0,\frac{1}{n}\right]$ and $m_1>\min\left\{2m_2 -(m_3+l), \max\{2m_2-1,\frac{n-2}{n}\}, m_2 -\frac{1}{n}\right\}$;
\item [] $\mathcal{A}_{ii}$) $\alpha \in \left(\frac{1}{n}, \frac{2}{n}\right)$ and 
$m_1>\min\left\{2m_2 +1-(m_3+l), \max\{2m_2,\frac{n-2}{n}\}, m_2 + \alpha-\frac{2}{n}\right\}$; 
\item [] $\mathcal{A}_{iii}$) $\alpha \in \left[\frac{2}{n},1\right]$ and
$m_1>\min\left\{2m_2 +1-(m_3+l), \max\{2m_2,\frac{n-2}{n}\}, m_2 + \frac{n\alpha-2}{n\alpha-1}\right\}$.
\end{itemize}
On the other hand, it is possible to check (see Example \ref{EsempioConfronto}) that $\mathcal{A}_i$) and $\mathcal{A}_{ii})$ are sharper than \ref{A3} and \ref{A4}, respectively. Conversely, $\mathcal{A}_{iii}$) and \ref{A5} are not directly comparable (the same applies for $\mathcal{A}_{iii}$) and \ref{Ap5}); see Figure \ref{fig:confronto}. 
Similarly, for the sake of completeness, the hypotheses
$\mathcal{A}_{iv})$ and $\mathcal{A}_{v})$ of \cite[Theorem 2.2]{FrassuLiViglialoro} become
\begin{itemize}
\item [] $\mathcal{A}_{iv}$) $\alpha \in \left(0,\frac{1}{n}\right]$ and  
$m_1>\min\left\{2m_2-(m_3+l), \max\{2m_2-1,\frac{n-2}{n}\}, 2m_2-\beta\right\}$,
\item [] $\mathcal{A}_{v}$) $\alpha \in \left(\frac{1}{n},1\right)$ and
$m_1>\min\left\{2m_2 +1-(m_3+l), \max\{2m_2,\frac{n-2}{n}\}, 2m_2 + 1 -\beta \right\}$,
\end{itemize}
and not even in this case any comparison indicating some sharper condition with respect to \ref{A1} and \ref{A2} can be observed.
\begin{figure}[h!]
 \centering
\begin{tikzpicture}
    \begin{axis}[hide axis,clip=false,xmin=-3,xmax=3,xlabel={X},ymin=-3,ymax=3,
        axis line style={->},
                 ytick=\empty,
                  xtick=\empty,
                   ztick=\empty,
      view={60}{45},
      samples=2,
      domain=-3:3,
      ]
    \addplot3[surf,fill=green,opacity=0.75]{3*x+2*y};                
    \addplot3[surf,fill=cyan,opacity=0.75]{y-2*x};   
    \end{axis}
  \coordinate (O) at (2,2,2);
  \draw[thick] (2,2,2) -- (2,2,2) node[below]{$O$};
  \draw[thick,->] (2,2,2) -- (7,2,4) node[below left=1]{$m_3$};
  \draw[thick,->] (2,2,2) -- (2,6.5,2) node[right=-1]{$\mathcal{A}_{iii}(m_2,m_3)$, $\mathcal{C}_{iii}(m_2,m_3)$};
  \draw[thick,->] (2,2,2) -- (2,2,7) node[right=2]{$m_2$};
    \end{tikzpicture}
     \caption{Comparing assumptions described in $\mathcal{A}_{iii}$) and \ref{A5}. For $n=3$, $\alpha=l=1$, we consider in terms of the variables $m_2$ and $m_3$ the functions 
$$\mathcal{A}_{iii}(m_2,m_3)=\min\left\{2m_2 +1-(m_3+l), \max\{2m_2,\frac{n-2}{n}\}, m_2 + \frac{n\alpha-2}{n\alpha-1}\right\}$$
and 
$$\mathcal{C}_{iii}(m_2,m_3)=\max\left\{m_2+l-\frac{2}{n},\min\left\{2m_3+1-(m_2+l),\max\left\{2m_3,\frac{n-2}{n}\right\}, m_3 + \frac{n\alpha-2}{n\alpha-1}\right\}\right\},$$
represented in green and cyan color, respectively. It is seen that in some regions of the plane $m_2Om_3$, there are values of $m_2$ and $m_3$ such that $\mathcal{A}_{iii}(m_2,m_3) > \mathcal{C}_{iii}(m_2,m_3)$ and other for which $\mathcal{A}_{iii}(m_2,m_3) < \mathcal{C}_{iii}(m_2,m_3)$.}
    \label{fig:confronto}
    \end{figure}
\begin{example}[Comparing assumptions \ref{A3} and  $\mathcal{A}_i)$]\label{EsempioConfronto}
Since the entire reasoning requires accurate computations, let us give here only some hints; in particular, let us indicate how to prove that assumption \ref{A3} is less mild than $\mathcal{A}_i)$. To this aim, let us use these abbreviations: $a=m_2+l-\frac{2}{n}$, $b=2m_3-(m_2+l)$, $c=\max\left\{2m_3-1,\frac{n-2}{n}\right\}$ and $d=m_3 - \frac{1}{n}$, $e=2m_2 -(m_3+l)$, $f=\max\{2m_2-1,\frac{n-2}{n}\}$ and $g= m_2 -\frac{1}{n}$. On the other hand, let $a=\max\{a,\min\{b,c,d\}\}$ and let us show that relation $a<\min\{e,f,g\}$ cannot be true; in this way nor $\min\{b,c,d\}<\min\{e,f,g\}$ may be satisfied. Indeed, if $a$ were smaller than $\min\{e,f,g\}$, then by imposing $a<e$ or $a<f$ or $a<g$ we would obtain these conditions contradicting our data: either $l<\frac{1}{n}$, or $l\geq 1$ but $m_2\in \emptyset.$ 
\end{example}
\end{remark}
\section{Local well posedness, boundedness criterion, main estimates and analysis of parameters}\label{SectionLocalInTime}
\textit{From now on we will tacitly assume that all the below constants $c_i$, $i=1,2,\ldots$ are positive.} 

For $\Omega$, $\chi,\xi,\delta$, $m_1,m_2,m_3$ and $f, g, h$ as above, $u, v, w \geq 0$ will indicate functions of $(x,t) \in \bar{\Omega}\times [0,\TM)$, for some finite $\TM$, classically solving problem \eqref{problem} when nonnegative initial data $(u_0,w_0)\in (W^{1,\infty}(\Omega))^2$ are provided.  In particular, $u$ satisfies 
\begin{equation}\label{massConservation}
\int_\Omega u(x, t)dx \leq m_0 \quad \textrm{for all }\, t \in (0,\TM), 
\end{equation}
whilst $w$ is such that 
\begin{equation}\label{wBounded}
0 \leq w\leq \lVert w_0\rVert_{L^\infty(\Omega)}\quad \textrm{in}\quad \Omega \times (0,\TM).
\end{equation}
Further, uniform boundedness of $(u,v,w)$ is ensured whenever (\textit{boundedness criterion}) $u\in L^\infty((0,\TM);L^p(\Omega))$, with $p>1$ arbitrarily large, for $m_1,m_2,m_3\in \R$, and with $p>\frac{n}{2}$, for $m_1=m_2=m_3=1$ (see, for instance, \cite{BellomoEtAl,TaoWinkParaPara}).

These statements are well-known in the mathematical context of chemotaxis models; as to the proof, when $h\equiv 0$ it verbatim follows from \cite[Lemmas 4.1 and 4.2]{FrassuCorViglialoro} and relation \eqref{massConservation} is the so-called  mass conservation property. Conversely, in the presence of the logistic term $h$ as in \eqref{h}, for $k_+=\max\{k,0\}$ an integration of the first equation in \eqref{problem} and an application of the H\"{o}lder inequality give 
\[
\frac{d}{dt} \int_\Omega u = \int_\Omega h(u) =k \int_\Omega u - \mu \int_\Omega u^{\beta} \leq k_+ \int_\Omega u - \frac{\mu}{|\Omega|^{\beta-1}} \left(\int_\Omega u\right)^{\beta}
\quad \textrm{for all }\, t \in (0,\TM),
\]
so that the uniform $L^1$-bound of $u$ is achieved by invoking an ODI-comparison argument.

In our computations, beyond the above positions, some uniform bounds of $\|w(\cdot,t)\|_{W^{1,s}(\Omega)}$ are required. In this sense, we mention the following lemma. (For the proof see \cite[Lemma 3.1]{FrassuLiViglialoro}.)
\begin{lemma}\label{LocalW}   
The $w$-component fulfills
\begin{equation}\label{Cg}
\int_\Omega |\nabla w(\cdot, t)|^s\leq \const{a} \quad \textrm{on } \,  (0,\TM)
\begin{cases}
\; \textrm{for all } s \in [1,\infty) & \textrm{if } \alpha \in \left(0, \frac{1}{n}\right],\\
\;  \textrm{for all } s \in \left[1, \frac{n}{n\alpha-1}\right) & \textrm{if } \alpha \in \left(\frac{1}{n},1\right].
\end{cases}  
\end{equation}
\end{lemma}
We will also make use of this technical result. 
\begin{lemma}\label{LemmaCoefficientAiAndExponents} 
Let $n\in \N$, with $n\geq 2$, $m_1>\frac{n-2}{n}$, $m_2,m_3\in \R$  and $\alpha \in (0,1]$. Then there is $s \in [1,\infty)$, such that for proper 
$p,r \in [1,\infty)$, $\theta$ and $\theta'$, $\mu$ and $\mu'$ conjugate exponents, we have that  
\begin{align}
a_1&= \frac{\frac{m_1+p-1}{2}\left(1-\frac{1}{(p+2m_3-m_1-1)\theta}\right)}{\frac{m_1+p-1}{2}+\frac{1}{n}-\frac{1}{2}},            
&  a_2&=\frac{r\left(\frac{1}{s}-\frac{1}{2\theta'}\right)}{\frac{r}{s}+\frac{1}{n}-\frac{1}{2}},  \nonumber \\ 
a_3 &= \frac{\frac{m_1+p-1}{2}\left(1-\frac{1}{2\alpha\mu}\right)}{\frac{m_1+p-1}{2}+\frac{1}{n}-\frac{1}{2}},  
&  a_4  &= \frac{r \left(\frac{1}{s}-\frac{1}{2(r-1)\mu'}\right)}{\frac{r}{s}+\frac{1}{n}-\frac{1}{2}},\nonumber\\ 
\kappa_1 &  =\frac{\frac{p}{2}\left(1- \frac{1}{p}\right)}{\frac{m_1+p-1}{2}+\frac{1}{n}-\frac{1}{2}} \nonumber, & \kappa_2 & =  \frac{r - \frac{1}{2}}{r+\frac{1}{n}-\frac{1}{2}}, 
\end{align}
belong to the interval $(0,1)$. If, indeed, 
\begin{equation}\label{Restrizionem1-m2-alphaPiccolo}
\alpha \in \left(0,\frac{1}{n}\right] \; \text{and}\; m_1>m_3-\frac{1}{n},
\end{equation}
\begin{equation}\label{Restrizionem1-m2-alphaGrande}
\alpha \in \left(\frac{1}{n},\frac{2}{n}\right) \; \text{and}\; m_1>m_3-\frac{2}{n}+\alpha,
\end{equation}
or
\begin{equation}\label{Restrizionem1-m2-alphaGrandeBis}
\alpha \in \left[\frac{2}{n},1\right] \; \text{and}\; m_1>m_3+\frac{n\alpha-2}{n\alpha-1},
\end{equation}
these further relations hold true: 
\begin{equation*}\label{MainInequalityExponents}
\beta_1 + \gamma_1 =\frac{p+2m_3-m_1-1}{m_1+p-1}a_1+\frac{1}{r}a_2\in (0,1) \;\textrm{ and }\;	\beta_2 + \gamma_2= \frac{2 \alpha }{m_1+p-1}a_3+\frac{r-1}{r}a_4 
\in (0,1).
\end{equation*}
\begin{proof}
All the details can be found in \cite[Lemma 3.2]{FrassuLiViglialoro}.
\end{proof}
\end{lemma}
\section{A priori estimates and proof of the Theorems}\label{EstimatesAndProofSection}
\subsection{The logistic case}\label{Log}
By relying on the globality criterion given in $\S$\ref{SectionLocalInTime}, let us dedicate to deriving the required uniform bound of $\int_\Omega u^p$. To this aim we define the functional $y(t):=\int_\Omega (u+1)^p + \int_\Omega |\nabla w|^{2p}$,  with $p>1$ properly large, and let us study the evolution in time of the functional $y(t)$ itself. In particular, we will soon see that during the analysis we will have to deal with terms of the form  $\int_{\Omega} (u+1)^{p+2m_3-m_1-1} \vert \nabla w\lvert^2$ and 
$\int_{\Omega} (u+1)^{2\alpha} \vert \nabla w\lvert^{2(p-1)}$. These integrals can be controlled  either invoking the Young inequality or the Gagliardo--Nirenberg one, and in this sense will be cornerstone the choice of $s$ introduced in Lemma \ref{LemmaCoefficientAiAndExponents}; specifically,  no restriction on the largeness of $s$ is possibly only when $\alpha \in \left(0,\frac{1}{n}\right]$.
\begin{lemma}\label{Estim_general_For_u^p_nablaw^2rLemmaLog}
If $m_1,m_2,m_3 \in \R$ comply with $m_1>2m_3-\beta$ and $m_2<\beta-l$ or $m_1>2m_3-(m_2+l)$  and $m_2<\beta-l$ or $m_1>\max\left\{2m_3-1,\frac{n-2}{n}\right\}$ 
and $m_2<\beta-l$ whenever $\alpha \in \left(0, \frac{1}{n}\right]$, or $m_1>2m_3-\beta+1$ and $m_2<\beta-l$ or $m_1>2m_3+1-(m_2+l)$ and $m_2<\beta-l$ or $m_1>\max\left\{2m_3,\frac{n-2}{n}\right\}$ and $m_2<\beta-l$ whenever $\alpha \in \left(\frac{1}{n},1\right)$, then there exists $p>1$ such that $(u,v,w)$ satisfies 
\begin{equation}\label{MainInequality_N}
\frac{d}{dt} \left(\int_\Omega (u+1)^p + \int_\Omega  |\nabla w|^{2p}\right) + c_{33} \int_\Omega |\nabla |\nabla w|^p|^2 + c_{34} \int_\Omega |\nabla (u+1)^{\frac{m_1+p-1}{2}}|^2 \leq c_{35} \quad \text{on } (0,\TM).
\end{equation}
\begin{proof}
Let $p>1$ and such that may be arbitrarily enlarged when necessary, and let us concentrate on the term $\frac{d}{dt} \int_\Omega (u+1)^p$. Testing the first equation of \eqref{problem} with $p(u+1)^{p-1}$ and exploiting assumption \eqref{h}, we conclude that
\begin{equation}\label{Estim1}
\begin{split}
&\frac{d}{dt} \int_\Omega (u+1)^p =\int_\Omega p(u+1)^{p-1}u_t \leq -p(p-1) \int_\Omega (u+1)^{p+m_1-3} |\nabla u|^2 
+p(p-1)\chi \int_\Omega u(u+1)^{m_2+p-3} \nabla u \cdot \nabla v \\
&-p(p-1)\xi \int_\Omega u(u+1)^{m_3+p-3} \nabla u \cdot \nabla w + pk \int_\Omega (u+1)^p - p \mu  \int_\Omega (u+1)^{p-1}u^{\beta} 
 \quad \text{for all } t \in (0,\TM).
\end{split}
\end{equation}
Now, setting $h_{p,m_2}(u)=p(p-1)\int_0^{u} \hat{u} (\hat{u}+1)^{m_2+p-3}\, d\hat{u}$, we have from $u<u+1$ that 
\begin{equation}\label{H}
\frac{p(p-1)}{p+m_2-1} u^{p+m_2-1} \leq h_{p,m_2}(u) \leq \frac{p(p-1)}{p+m_2-1} [(u+1)^{p+m_2-1}-1] \quad \textrm{in} \quad \Omega \times (0,\TM). 
\end{equation}
Successively, we focus on the second integral of the rhs of \eqref{Estim1}; by taking into account the growth of $g$ specified in \eqref{f}, the inequality $(A+B)^p \leq 2^{p-1} (A^p+B^p)$, with $A,B \geq 0$ and $p>1$ (it being understood its employment in the next lines), and bound \eqref{H}, we arrive for $p$ sufficiently large at
\begin{equation}\label{Estim2}
\begin{split}
&p(p-1)\chi \int_{\Omega} u(u+1)^{m_2+p-3} \nabla u \cdot \nabla v = -\chi \int_\Omega h_{p,m_2}(u) \Delta v \\
&\leq -\frac{p(p-1)\chi \delta}{p+m_2-1} \int_\Omega u^{p+m_2-1}v + \frac{p(p-1)\chi}{p+m_2-1} \int_\Omega [(u+1)^{p+m_2-1}-1]g(u)\\
&\leq  -\frac{p(p-1)\chi \delta}{p+m_2-1} \int_\Omega \left[\frac{(u+1)^{p+m_2-1}}{2^{p+m_2-2}}-1\right]v + \frac{p(p-1)\chi}{p+m_2-1} \int_\Omega (u+1)^{p+m_2-1}g(u) - \frac{p(p-1)\chi}{p+m_2-1} \int_\Omega g(u)\\
&\leq -\frac{p(p-1)\chi \delta}{2^{p+m_2-2}(p+m_2-1)} \int_\Omega (u+1)^{p+m_2-1} v + \frac{p(p-1)\chi \gamma}{p+m_2-1} \int_\Omega (u+1)^l\\
&+\frac{p(p-1)\chi \gamma}{p+m_2-1} \int_\Omega (u+1)^{p+m_2+l-1} -\frac{p(p-1)\chi \gamma}{p+m_2-1} \int_\Omega u^l\\
& \leq -\frac{p(p-1)\chi \delta}{2^{p+m_2-2}(p+m_2-1)} \int_\Omega (u+1)^{p+m_2-1} v + \frac{p(p-1)\chi \gamma}{p+m_2-1} \int_\Omega (u+1)^{p+m_2+l-1}\\
&-\frac{p(p-1)\chi \gamma}{p+m_2-1} \left(\frac{1}{2^{l-1}}-1\right) \int_\Omega (u+1)^l + \frac{p(p-1)\chi \gamma |\Omega|}{p+m_2-1}\\
& \leq \frac{p(p-1)\chi \gamma}{p+m_2-1} \int_\Omega (u+1)^{p+m_2+l-1}
-\frac{p(p-1)\chi \gamma}{p+m_2-1} \left(\frac{1}{2^{l-1}}-1\right) \int_\Omega (u+1)^l + \frac{p(p-1)\chi \gamma |\Omega|}{p+m_2-1} \quad \textrm{on } (0,\TM),
\end{split}
\end{equation}
where we made use of the identity $\delta \int_\Omega v = \int_\Omega g(u)$, naturally coming  by integrating the second equation of \eqref{problem}.

Regarding the second integral of the rhs of \eqref{Estim2}, we note that if $l=1$ such term vanishes, while for $l>1$ and $p$ sufficiently large it can be treated by Young's inequality, entailing 
\begin{equation}\label{L0}
-\frac{p(p-1)\chi \gamma}{p+m_2-1} \left(\frac{1}{2^{l-1}}-1\right) \int_\Omega (u+1)^l \leq \epsilon \int_\Omega (u+1)^{p-1+\beta} + \const{bl} \quad \textrm{on } (0,\TM).
\end{equation}
Applying again Young's inequality to the third integral of the rhs of \eqref{Estim1} yields for $\epsilon>0$ 
\begin{equation}\label{Young}
-p(p-1)\xi \int_\Omega u(u+1)^{m_3+p-3} \nabla u \cdot \nabla w \leq \epsilon \int_\Omega (u+1)^{p+m_1-3}|\nabla u|^2 
+ \const{c} \int_\Omega (u+1)^{p+2m_3-m_1-1}|\nabla w|^2 \quad \textrm{on } (0,\TM).
\end{equation}
From $-u^{\beta} \leq -\frac{1}{2^{\beta-1}} (u+1)^{\beta}+1$, the last integral in \eqref{Estim1} is transformed into
\begin{equation}\label{beta}
- p \mu  \int_\Omega (u+1)^{p-1} u^{\beta} \leq -\frac{p \mu}{2^{\beta-1}} \int_\Omega (u+1)^{p-1+\beta} + p \mu \int_\Omega (u+1)^{p-1} \quad \text{for all } t \in (0,\TM).
\end{equation}
Henceforth, a further application of the Young inequality yields for all $t \in (0,\TM)$
\begin{equation} \label{k}
pk \int_\Omega (u+1)^p \leq pk_+ \int_\Omega (u+1)^p \leq \epsilon \int_\Omega (u+1)^{p-1+\beta} + \const{d} \quad \text{and} \quad  
p \mu \int_\Omega (u+1)^{p-1} \leq \epsilon \int_\Omega (u+1)^{p-1+\beta} + \const{e}.
\end{equation}
Let us turn our attention to the term $\frac{d}{dt} \int_\Omega  |\nabla w|^{2p}$ of the functional $y(t)$; herein we avoid giving all the details, which can be found in \cite[Lemma 5.3]{FrassuCorViglialoro}, and we only say that after some computations one can deduce that  
\begin{equation}\label{Estim_gradW}
\frac{d}{dt} \int_\Omega  |\nabla w|^{2p}+ p \int_\Omega |\nabla w|^{2p-2} |D^2w|^2 \leq \const{f} \int_\Omega u^{2\alpha} |\nabla w|^{2p-2} + \const{g} \quad \textrm{on } \;  (0,\TM).
\end{equation}
By putting together  estimates \eqref{Estim1} and \eqref{Estim_gradW}, and in view of the gained bounds \eqref{Estim2}, \eqref{L0}, \eqref{Young}, \eqref{beta} and \eqref{k}, we can write 
\begin{equation}\label{EstimSum}
\begin{split}
&\frac{d}{dt} \left(\int_\Omega (u+1)^p +  \int_\Omega  |\nabla w|^{2p}\right) + p \int_\Omega |\nabla w|^{2p-2} |D^2w|^2 
\leq (\epsilon-p(p-1)) \int_\Omega (u+1)^{p+m_1-3} |\nabla u|^2\\
&+ \left(3\epsilon-\frac{p\mu}{2^{\beta-1}}\right) \int_\Omega (u+1)^{p-1+\beta}+ \const{c} \int_\Omega (u+1)^{p+2m_3-m_1-1}|\nabla w|^2\\
&+ \const{f} \int_\Omega u^{2\alpha} |\nabla w|^{2p-2} + \frac{p(p-1)\chi \gamma}{p+m_2-1} \int_\Omega (u+1)^{p+m_2+l-1} +\const{h} \quad \text{for all } t \in (0,\TM).
\end{split}
\end{equation}
Finally, in order to treat the integrals with $\nabla w$ on the rhs of \eqref{EstimSum}, let us  distinguish two cases, each one correlated to a specific range of $\alpha$.
\begin{itemize}
\item \textbf{Case $1$}: $\alpha \in \left(0, \frac{1}{n}\right]$ and $m_1>2m_3-\beta$ and $m_2<\beta-l$ or $m_1>2m_3-(m_2+l)$ and $m_2<\beta-l$ or $m_1>\max\left\{2m_3-1,\frac{n-2}{n}\right\}$ and $m_2<\beta-l$. With $s$ arbitrarily large, an application of the Young inequality and bound \eqref{Cg} provide on $(0,\TM)$
\begin{equation}\label{YoungS}
\begin{split}
\const{c} \int_\Omega (u+1)^{p+2m_3-m_1-1}|\nabla w|^2 &\leq \int_\Omega |\nabla w|^{s} + \const{i} \int_\Omega (u+1)^{\frac{(p+2m_3-m_1-1)s}{s-2}}\leq \const{i} \int_\Omega (u+1)^{\frac{(p+2m_3-m_1-1)s}{s-2}} + \const{j}.
\end{split}
\end{equation}
Due to $m_1>2m_3-\beta$, we have $\frac{(p+2m_3-m_1-1)s}{s-2}<p-1+\beta$, so that Young's inequality yields 
\begin{equation}\label{YoungS1}
\const{i} \int_\Omega (u+1)^{\frac{(p+2m_3-m_1-1)s}{s-2}} \leq \epsilon \int_\Omega (u+1)^{p-1+\beta} + \const{l}  \quad \textrm{for all }\, t \in (0,\TM).
\end{equation}
On the other hand, from $m_1>2m_3-(m_2+l)$, it is obtained $\frac{(p+2m_3-m_1-1)s}{s-2}<m_2+p+l-1$, and another application of the Young inequality implies
\begin{equation}\label{YoungS2}
\const{i} \int_\Omega (u+1)^{\frac{(p+2m_3-m_1-1)s}{s-2}} \leq   \int_\Omega (u+1)^{m_2+p+l-1} + \const{m}  \quad \textrm{on }\, (0,\TM).
\end{equation}
Similarly, from $m_1>2m_3-1$ it holds that $\frac{(p+2m_3-m_1-1)s}{s-2}<p$, and subsequently
\begin{equation}\label{YoungS3}
\const{i} \int_\Omega (u+1)^{\frac{(p+2m_3-m_1-1)s}{s-2}} \leq  \int_\Omega (u+1)^p + \const{n}  \quad \textrm{for all }\, t \in (0,\TM).
\end{equation}
At this point, by the Gagliardo--Nirenberg inequality and the uniform-in-time $L^1$ bound \eqref{massConservation} we have 
\begin{equation*}
\begin{split}
\int_{\Omega} (u+1)^p&= \|(u+1)^{\frac{m_1+p-1}{2}}\|_{L^{\frac{2p}{m_1+p-1}}(\Omega)}^{\frac{2p}{m_1+p-1}}\\ 
&\leq \const{o} \|\nabla (u+1)^{\frac{m_1+p-1}{2}}\|_{L^2(\Omega)}^{\frac{2p}{m_1+p-1}\theta_1} \|(u+1)^{\frac{m_1+p-1}{2}}\|_{L^{\frac{2}{m_1+p-1}}(\Omega)}^{\frac{2p}{m_1+p-1}(1-\theta_1)} + \const{o} \|(u+1)^{\frac{m_1+p-1}{2}}\|_{L^{\frac{2}{m_1+p-1}}(\Omega)}^{\frac{2p}{m_1+p-1}}\\
& \leq \const{p} \Big(\int_\Omega |\nabla (u+1)^\frac{m_1+p-1}{2}|^2\Big)^{\kappa_1}+ \const{p} \quad \text{ for all } t \in(0,\TM),
\end{split}
\end{equation*}
where
\begin{equation*}
\theta_1=\frac{\frac{n(m_1+p-1)}{2}\left(1-\frac{1}{p}\right)}{1-\frac{n}{2}+\frac{n(m_1+p-1)}{2}}\in (0,1).
\end{equation*}
Since $\kappa_1 \in (0,1)$ from Lemma \ref{LemmaCoefficientAiAndExponents}, the Young inequality produces 
\begin{equation}\label{YoungS4}
\int_\Omega (u+1)^p \leq \epsilon \int_\Omega |\nabla (u+1)^\frac{m_1+p-1}{2}|^2  + \const{q} \quad \text{on } (0,\TM).
\end{equation}
Finally, by exploiting twice the Young inequality (for $p,s$ large) and bound \eqref{Cg}, we obtain the following estimate: 
\begin{equation}\label{Estimat_nablaw^2r+2}
\begin{split}
&\const{f} \int_\Omega u^{2\alpha} |\nabla w|^{2p-2} \leq \epsilon \int_\Omega (u+1)^{p-1+\beta} + \const{r} \int_\Omega |\nabla w|^{\frac{2(p-1)(p-1+\beta)}{p-1+\beta-2\alpha}}\\
& \leq \epsilon \int_\Omega (u+1)^{p-1+\beta} + \int_\Omega |\nabla w|^s + \const{s} 
\leq \epsilon \int_\Omega (u+1)^{p-1+\beta} + \const{t} \quad \textrm{for all} \quad t \in (0,\TM).
\end{split}
\end{equation}
Alternatively, the previous estimate can be also reformulated as
\begin{equation}\label{Estimat_nablaW^2r+2}
\begin{split}
&\const{f} \int_\Omega u^{2\alpha} |\nabla w|^{2p-2} \leq  \int_\Omega (u+1)^{m_2+p+l-1} + \const{u} \int_\Omega |\nabla w|^{\frac{2(p-1)(m_2+p+l-1)}{m_2+p+l-1-2\alpha}}\\
& \leq  \int_\Omega (u+1)^{m_2+p+l-1} + \int_\Omega |\nabla w|^s + \const{bo} 
\leq  \int_\Omega (u+1)^{m_2+p+l-1} + \const{v} \quad \textrm{on} \quad (0,\TM).
\end{split}
\end{equation}
\item  
\textbf{Case $2$}: $\alpha \in \left(\frac{1}{n},1\right)$ and $m_1>2m_3-\beta+1$ and $m_2<\beta-l$ or $m_1>2m_3+1-(m_2+l)$ and $m_2<\beta-l$ or $m_1>\max\left\{2m_3,\frac{n-2}{n}\right\}$ and $m_2<\beta-l$. In this case, $s$ cannot arbitrarily increase (it belongs to a bounded interval), so we have to manipulate the integral \eqref{YoungS} and the last one in \eqref{Estim_gradW} in an alternative  way. To be precise, a double application of the Young inequality gives  
\begin{equation}\label{YoungS5}
\begin{split}
&\const{c} \int_\Omega (u+1)^{p+2m_3-m_1-1}|\nabla w|^2 \leq \epsilon \int_\Omega |\nabla w|^{2(p+1)} + \const{z} \int_\Omega (u+1)^{\frac{(p+2m_3-m_1-1)(p+1)}{p}}\\
& \leq \epsilon  \int_\Omega |\nabla w|^{2(p+1)} + \epsilon \int_\Omega (u+1)^{p-1+\beta} + \const{Miso} \quad \textrm{for all} \quad t \in (0,\TM).
\end{split}
\end{equation}
From $m_1>2m_3+1-(m_2+l)$ we also have $\frac{(p+2m_3-m_1-1)(p+1)}{p}<m_2+p+l-1$, and Young's inequality entails with the property that 
\begin{equation}\label{YoungS6}
\const{z} \int_\Omega (u+1)^{\frac{(p+2m_3-m_1-1)(p+1)}{p}} \leq  \int_\Omega (u+1)^{m_2+p+l-1} + \const{ab} \quad \textrm{on} \quad (0,\TM).
\end{equation}
But, in view of the fact that $m_1>\max\{2m_3,\frac{n-2}{n}\}$ implies $\frac{(p+2m_3-m_1-1)(p+1)}{p}<p$ one can estimate the previous integral by applying the Young and Gagliardo--Nirenberg inequalities; we have, as already done above,
\begin{equation}\label{YoungS7}
\const{z} \int_\Omega (u+1)^{\frac{(p+2m_3-m_1-1)(p+1)}{p}} \leq  \int_\Omega (u+1)^p + \const{ac}
\leq \epsilon \int_\Omega |\nabla (u+1)^\frac{m_1+p-1}{2}|^2 + \const{ad} \quad \textrm{for all} \quad t \in (0,\TM).
\end{equation}
Similarly, also for the integral in the rhs of \eqref{Estim_gradW}, we can obtain two possible estimates; indeed, by exploiting Young's inequality and $\alpha<1$ we have  
\begin{equation}\label{Estimat_nablaw^2r+2_1}
\begin{split}
&\const{f} \int_\Omega u^{2\alpha} |\nabla w|^{2p-2} \leq \epsilon \int_\Omega (u+1)^{p-1+\beta} + \const{ae} \int_\Omega |\nabla w|^{\frac{2(p-1)(p-1+\beta)}{p-1+\beta-2\alpha}}\\
& \leq \epsilon \int_\Omega (u+1)^{p-1+\beta} + \epsilon \int_\Omega |\nabla w|^{2(p+1)} + \const{af} \quad \textrm{for all} \quad t \in (0,\TM),
\end{split}
\end{equation}
or 
\begin{equation}\label{Estimat_nablaW^2r+2_1}
\begin{split}
&\const{f} \int_\Omega u^{2\alpha} |\nabla w|^{2p-2} \leq  \int_\Omega (u+1)^{m_2+p+l-1} + \const{ag} \int_\Omega |\nabla w|^{\frac{2(p-1)(m_2+p+l-1)}{m_2+p+l-1-2\alpha}}\\
& \leq  \int_\Omega (u+1)^{m_2+p+l-1}  + \epsilon \int_\Omega |\nabla w|^{2(p+1)} + \const{ah} \quad \textrm{on} \quad (0,\TM).
\end{split}
\end{equation}
In turn, by invoking \cite[(9), page 6]{FrassuCorViglialoro} and bound \eqref{wBounded} for the term involving $\int_\Omega |\nabla w|^{2(p+1)}$, we get
\begin{equation}\label{GradSfin}
\int_\Omega |\nabla w|^{2(p+1)} \leq 2 (4p^2+n) \|w_0\|_{L^{\infty}(\Omega)}^2 \int_\Omega |\nabla w|^{2p-2} |D^2 w|^2 \quad \textrm{on } (0,\TM).
\end{equation}
Finally, as to the last integral of the rhs of \eqref{EstimSum}, by exploiting the Young inequality and $m_2<\beta-l$, we get for any $\tilde{c}>0$ 
\begin{equation}\label{Young0}
\tilde{c} \int_\Omega (u+1)^{p+m_2+l-1} \leq \epsilon \int_\Omega (u+1)^{p-1+\beta} + \const{b} \quad \textrm{for all } t \in (0,\TM).
\end{equation}
As a consequence, by rephrasing estimate \eqref{EstimSum} in terms of bounds \eqref{YoungS}--\eqref{Young0}, a proper choice of the free parameter $\epsilon>0$ gives 
\[
\frac{d}{dt} \left(\int_\Omega (u+1)^p + \int_\Omega  |\nabla w|^{2p}\right) + \const{ai} \int_\Omega |\nabla |\nabla w|^p|^2 + \const{aj} \int_\Omega |\nabla (u+1)^{\frac{m_1+p-1}{2}}|^2 \leq \const{al} \quad \text{on } (0,\TM),
\]
where we also exploited that 
\begin{equation*}
\int_\Omega (u+1)^{p+m_1-3} |\nabla u|^2 = \frac{4}{(m_1+p-1)^2} \int_\Omega |\nabla (u+1)^{\frac{m_1+p-1}{2}}|^2 \quad  \text{on } (0,\TM),
\end{equation*}
and (see \cite[page 17]{FrassuCorViglialoro})
\begin{equation*}
\vert \nabla \lvert \nabla w\rvert^p\rvert^2=\frac{p^2}{4}\lvert \nabla w \rvert^{2p-4}\vert \nabla \lvert \nabla w\rvert^2\rvert^2=p^2\lvert \nabla w \rvert^{2p-4}\lvert D^2w \nabla w \rvert^2\leq p^2|\nabla w|^{2p-2} |D^2w|^2.
\end{equation*}
Let us give a final comment: the assumptions of this lemma can be presented in this more compact formulation:
\begin{equation*}
m_2<\beta-l \textrm{ and }
\begin{cases}
 m_1>\min\left\{2m_3-\beta, 2m_3-(m_2+l),\max\left\{2m_3-1,\frac{n-2}{n}\right\}\right\} & \textrm{ for } \alpha \in \left(0, \frac{1}{n}\right],\\
 m_1>\min\left\{2m_3-\beta+1, 2m_3+1-(m_2+l),\max\left\{2m_3,\frac{n-2}{n}\right\}\right\} & \textrm{ for } \alpha \in \left(\frac{1}{n}, 1\right).
\end{cases}
\end{equation*}
In particular, since for $m_2<\beta-l$ it is seen that $2m_3-\beta< 2m_3-(m_2+l)$ and  $2m_3-\beta+1< 2m_3+1-(m_2+l)$, eventually we have that these assumptions exactly read as in Theorem \ref{MainTheorem1}.
\end{itemize}
Let us conclude by considering the case $\alpha=1$. It is seen that  both bounds \eqref{Estimat_nablaw^2r+2_1} and \eqref{Estimat_nablaW^2r+2_1} can be identically derived by means of the Young inequality, whenever $\beta>2$ (for the first) and $\beta>2$ and $m_2>2-l$ (for the second; recall $m_2<\beta-l$).
\end{proof}
\end{lemma}
\subsection{The non-logistic case}\label{NonLog}
\begin{lemma}\label{Estim_general_For_u^p_nablaw^2rLemma}
If $m_1,m_2,m_3 \in \R$ comply with $m_1>2m_3-(m_2+l)$ and $m_1>m_2+l-\frac{2}{n}$ or $m_1>\max\left\{2m_3-1,\frac{n-2}{n}\right\}$ and $m_1>m_2+l-\frac{2}{n}$ whenever $\alpha \in \left(0, \frac{1}{n}\right]$, or $m_1>2m_3+1-(m_2+l)$ and $m_1>m_2+l-\frac{2}{n}$ or $m_1>\max\left\{2m_3,\frac{n-2}{n}\right\}$ and $m_1>m_2+l-\frac{2}{n}$ whenever $\alpha \in \left(\frac{1}{n},1\right)$, then there exists $p>1$ such that $(u,v,w)$ satisfies a similar inequality as in \eqref{MainInequality_N}.
\begin{proof}
By reasoning as in the proof of Lemma \ref{Estim_general_For_u^p_nablaw^2rLemmaLog}, estimate \eqref{EstimSum} becomes 
\begin{equation}\label{Estim1N}
\begin{split}
&\frac{d}{dt} \left(\int_\Omega (u+1)^p +  \int_\Omega  |\nabla w|^{2p}\right) + p \int_\Omega |\nabla w|^{2p-2} |D^2w|^2 
\leq (\epsilon-p(p-1)) \int_\Omega (u+1)^{p+m_1-3} |\nabla u|^2\\ 
&+ \frac{p(p-1)\chi \gamma}{p+m_2-1} \int_\Omega (u+1)^{p+m_2+l-1} -\frac{p(p-1)\chi \gamma}{p+m_2-1} \left(\frac{1}{2^{l-1}}-1\right) \int_\Omega (u+1)^l \\
&+ \const{c} \int_\Omega (u+1)^{p+2m_3-m_1-1}|\nabla w|^2
+ \const{f} \int_\Omega u^{2\alpha} |\nabla w|^{2p-2}+ \const{am} \quad \text{for all } 
t \in (0,\TM).
\end{split}
\end{equation}
As to the third integral of the rhs of \eqref{Estim1N}, as we have seen before, if $l=1$ such term vanishes, while for $l>1$ Young's inequality provides 
\begin{equation}\label{L} 
-\frac{p(p-1)\chi \gamma}{p+m_2-1} \left(\frac{1}{2^{l-1}}-1\right) \int_\Omega (u+1)^l
\leq \epsilon \int_\Omega (u+1)^{p+m_2+l-1} + \const{bbl} \quad \textrm{on } (0,\TM).
\end{equation}
Successively, thanks to the Gagliardo--Nirenberg inequality and \eqref{massConservation} (for $h\equiv 0$ this is precisely the mass conservation property for $u$), the estimate of the second integral of the rhs in \eqref{Estim1N} also reads for any $\hat{c}>0$
\begin{equation*}
\begin{split}
&\hat{c} \int_\Omega (u+1)^{p+m_2+l-1}= \hat{c} \|(u+1)^{\frac{p+m_1-1}{2}}\|_{L^{\frac{2(p+m_2+l-1)}{p+m_1-1}}(\Omega)}^{\frac{2(p+m_2+l-1)}{p+m_1-1}}\\
& \leq \const{an} \|\nabla (u+1)^{\frac{p+m_1-1}{2}}\|_{L^2(\Omega)}^{\frac{2(p+m_2+l-1)}{p+m_1-1}\theta_2} 
\|(u+1)^{\frac{p+m_1-1}{2}}\|_{L^{\frac{2}{p+m_1-1}}(\Omega)}^{\frac{2(p+m_2+l-1)}{p+m_1-1}(1-\theta_2)} 
+ \const{an} \|(u+1)^{\frac{p+m_1-1}{2}}\|_{L^{\frac{2}{p+m_1-1}}(\Omega)}^{\frac{2(p+m_2+l-1)}{p+m_1-1}}\\
&\leq  \const{ao} \left(\int_\Omega |\nabla (u+1)^{\frac{p+m_1-1}{2}}|^2\right)^{\frac{(p+m_2+l-1)}{p+m_1-1}\theta_2} + \const{ao} \quad \textrm{on } (0,\TM),
\end{split}
\end{equation*}
where
\[
\theta_2=\frac{\frac{p+m_1-1}{2} - \frac{p+m_1-1}{2(p+m_2+l-1)}}{\frac{p+m_1-1}{2}+\frac{1}{n}-\frac{1}{2}} \in (0,1).
\]
Since $\frac{(p+m_2+l-1)}{p+m_1-1}\theta_2<1$ for $m_1>m_2 + l-\frac{2}{n}$, an application of the Young inequality gives 
\begin{equation}\label{GN1_N}
\hat{c} \int_\Omega (u+1)^{p+m_2+l-1} \leq \epsilon \int_\Omega |\nabla (u+1)^{\frac{p+m_1-1}{2}}|^2 + \const{ap} \quad \textrm{for all } t \in (0,\TM).
\end{equation}
Now, as before, we have to distinguish
\begin{itemize}
\item \textbf{Case $1$}: $\alpha \in \left(0, \frac{1}{n}\right]$ and $m_1>2m_3-(m_2+l)$ and $m_1>m_2+l-\frac{2}{n}$ or $m_1>\max\left\{2m_3-1,\frac{n-2}{n}\right\}$ and $m_1>m_2+l-\frac{2}{n}$,
\item \textbf{Case $2$}: $\alpha \in \left(\frac{1}{n},1\right)$ and $m_1>2m_3+1-(m_2+l)$ and $m_1>m_2+l-\frac{2}{n}$ or $m_1>\max\left\{2m_3,\frac{n-2}{n}\right\}$ and $m_1>m_2+l-\frac{2}{n}$,
\end{itemize} 
and both have been already analyzed in the proof of Lemma \ref{Estim_general_For_u^p_nablaw^2rLemmaLog};  therefore, by following the same arguments we have the claim. 

Also in this case, we point out that Lemma \ref{Estim_general_For_u^p_nablaw^2rLemma} and in particular inequality \eqref{Estimat_nablaW^2r+2_1} is satisfied also in the limit case $\alpha=1$, provided that the condition $m_2>2-l$ is required.
\end{proof}
\end{lemma}
In addition, in order to move toward the proofs of our claims, by employing the Gagliardo--Nirenberg inequality, we as well establish this result.
\begin{lemma}\label{Met_GN}
If $m_1,m_2, m_3\in \R$ and $\alpha>0$ are taken accordingly to \eqref{Restrizionem1-m2-alphaPiccolo}, \eqref{Restrizionem1-m2-alphaGrande}, 
\eqref{Restrizionem1-m2-alphaGrandeBis} and to $m_1>m_2+l-\frac{2}{n}$, then there exist $p, r>1$ such that $(u,v,w)$ satisfies a similar inequality as in \eqref{MainInequality_N}.
\begin{proof}
Herein we have to consider a functional slightly modified with respect $y(t)$; more precisely, it reads $\int_\Omega (u+1)^p + \int_\Omega |\nabla w|^{2r}$, with $p, r>1$ properly large, and generally with $p\neq r$. In this way, estimate \eqref{Estim1N} appears as 
\begin{equation*}
\begin{split}
&\frac{d}{dt} \left(\int_\Omega (u+1)^p +  \int_\Omega  |\nabla w|^{2r}\right) + r \int_\Omega |\nabla w|^{2r-2} |D^2w|^2 
\leq (\epsilon-p(p-1)) \int_\Omega (u+1)^{p+m_1-3} |\nabla u|^2\\ 
&+ \frac{p(p-1)\chi \gamma}{p+m_2-1} \int_\Omega (u+1)^{p+m_2+l-1} -\frac{p(p-1)\chi \gamma}{p+m_2-1} \left(\frac{1}{2^{l-1}}-1\right) \int_\Omega (u+1)^l \\
&+ \const{c} \int_\Omega (u+1)^{p+2m_3-m_1-1}|\nabla w|^2
+ \const{f} \int_\Omega u^{2\alpha} |\nabla w|^{2r-2}+ \const{am} \quad \text{for all } 
t \in (0,\TM).
\end{split}
\end{equation*}
As to the terms involving $\nabla w$, they are estimated exactly as done in the proof of \cite[Lemma 4.2]{FrassuLiViglialoro}; in particular, the coefficients $a_1,a_2,a_3,a_4,\beta_1,\gamma_1,\beta_2$ and $\gamma_2$ defined in Lemma \ref{LemmaCoefficientAiAndExponents} take part in the related  logical steps. On the other hand, and conversely to \cite[Lemma 4.2]{FrassuLiViglialoro}, the additional two terms proportional to $\int_\Omega (u+1)^l$ and $\int_\Omega (u+1)^{p+m_2+l-1}$ are controlled respectively by means of bounds \eqref{L} and \eqref{GN1_N} once, in this order, $p$ is taken arbitrarily large and $m_1>m_2+l-\frac{2}{n}$ is considered. 
\end{proof}
\end{lemma}
As a by-product of what now obtained we are in a position to conclude.
\subsection{Proof of Theorems \ref{MainTheorem} and \ref{MainTheorem1}}
\begin{proof}
Let $(u_0,w_0) \in (W^{1,\infty}(\Omega))^2$ with $u_0, w_0 \geq 0$ on $\bar{\Omega}$. For $f$, $g$ as in \eqref{f} and $h$ as in \eqref{h}, let $\alpha >0$ and let $m_1,m_2,m_3 \in \R$ comply with \ref{A1} and \ref{A2}, respectively, \ref{A3}, \ref{A4} and \ref{A5}. Then, we refer to Lemma
\ref{Estim_general_For_u^p_nablaw^2rLemmaLog}, and Lemmas \ref{Estim_general_For_u^p_nablaw^2rLemma}, \ref{Met_GN}, in this order, and we obtain for some 
$C_1,C_2,C_3>0$
\begin{equation}\label{Estim_general_For_y_2}
y'(t) + C_1 \int_\Omega |\nabla (u+1)^{\frac{m_1+p-1}{2}}|^2 + C_2 \int_\Omega |\nabla |\nabla w|^p|^2 \leq C_3
\quad \text{ on } (0, \TM).
\end{equation}
Successively, the Gagliardo--Nirenberg inequality again implies this estimate
\begin{equation*}\label{G_N2}
\int_\Omega (u+1)^p \leq \const{aq} \Big(\int_\Omega |\nabla (u+1)^\frac{m_1+p-1}{2}|^2\Big)^{\kappa_1}+ \const{aq}  \quad \text{for all } t \in (0,\TM),
\end{equation*}
(as already done in \eqref{YoungS4}), and also the following one
\begin{equation*}
\int_\Omega \lvert \nabla w\rvert^{2p}=\lvert \lvert \lvert \nabla w\rvert^p\lvert \lvert_{L^2(\Omega)}^2 
\leq \const{ao} \lvert \lvert\nabla  \lvert \nabla w \rvert^p\rvert \lvert_{L^2(\Omega)}^{2\kappa_2} \lvert \lvert\lvert \nabla w \rvert^p\lvert \lvert_{L^\frac{1}{p}(\Omega)}^{2(1-\kappa_2)} + \const{ao} \lvert \lvert \lvert \nabla w \rvert^p\lvert \lvert^2_{L^\frac{1}{p}(\Omega)}\quad \textrm{on } (0,\TM),
\end{equation*}
with $\kappa_2$ taken from  Lemma \ref{LemmaCoefficientAiAndExponents}.  Subsequently, the uniform-in-time $L^s$-bound of $\nabla w$ in \eqref{Cg} infers 
\begin{equation*}\label{Estim_Nabla nabla v^p^2}
\int_\Omega \lvert \nabla w\rvert^{2p}\leq \const{ar} \Big(\int_\Omega \lvert \nabla \lvert \nabla w \rvert^p\rvert^2\Big)^{\kappa_2}+\const{ar} \quad \text{for all } t \in (0,\TM).
\end{equation*}
In the same spirit of \cite[Lemma 5.4]{FrassuCorViglialoro}, by using the estimates involving $\int_\Omega (u+1)^p$ and $\int_\Omega \rvert\nabla w\lvert^{2p}$,
relation \eqref{Estim_general_For_y_2} provides $\kappa=\min\{\frac{1}{\kappa_1},\frac{1}{\kappa_2}\}$ such that 
\begin{equation*}\label{MainInitialProblemWithM}
\begin{cases}
y'(t)\leq \const{as}-\const{at} y^{\kappa}(t)\quad \textrm{for all } t \in (0,\TM),\\
y(0)=\int_\Omega (u_0+1)^p+ \int_\Omega |\nabla w_0|^{2p}, 
\end{cases}
\end{equation*}
so giving some $L>0$ such that, through comparison arguments, $\int_\Omega (u+1)^p \leq L$ for all $t\in (0,\TM)$. Finally, the gained inclusion $u \in L^\infty((0,\TM);L^p(\Omega))$ and the boundedness criterion in 
$\S$\ref{SectionLocalInTime} give the statement. 
\end{proof}
\subsection{Proof of Corollary \ref{MainCorollary1}}
\begin{proof}
The first claim of Corollary \ref{MainCorollary1} is a direct consequence of Theorem \ref{MainTheorem1}, simply by substituting $m_1=m_2=m_3=1$ in \ref{A1} and \ref{A2}. For the second part, for $m_1=m_2=m_3=1$ and $\beta=1+l$ estimate \eqref{EstimSum} becomes, by  explicitly computing the value of $c_3$,
\begin{equation}\label{Est1L}
\begin{split}
&\frac{d}{dt} \left(\int_\Omega (u+1)^p + \int_\Omega |\nabla w|^{2p}\right) + p \int_\Omega |\nabla w|^{2p-2} |D^2 w|^2 \leq 
-\frac{p(p-1)}{2} \int_\Omega (u+1)^{p-2} |\nabla u|^2 + \frac{\xi^2 p(p-1)}{2} \int_\Omega (u+1)^p |\nabla w|^2\\ 
&+ \const{la} \int_\Omega u^{2\alpha} |\nabla w|^{2p-2} + \left((p-1)\chi \gamma +3\epsilon-\frac{\mu p}{2^l}\right) \int_\Omega (u+1)^{p+l}+ \const{lb} \quad \text{for all } 
t \in (0,\TM).
\end{split}
\end{equation}
Thanks to the Young inequality, we obtain that the second and the third integral of \eqref{Est1L} satisfy for $\epsilon>0$
\begin{equation}\label{Lin1}
\frac{\xi^2 p(p-1)}{2} \int_\Omega (u+1)^p |\nabla w|^2 \leq \epsilon \int_\Omega (u+1)^{p+l} + \const{lc} \int_\Omega |\nabla w|^{\frac{2(p+l)}{l}} \quad \text{on } (0,\TM),
\end{equation}
and
\begin{equation}\label{Lin2}
\const{la} \int_\Omega u^{2\alpha} |\nabla w|^{2p-2} \leq \epsilon \int_\Omega (u+1)^{p+l} + \const{ld} \int_\Omega |\nabla w|^{\frac{2(p-1)(p+l)}{p+l-2\alpha}} \quad \text{for all } t \in (0,\TM).
\end{equation}
At this point, in order to treat the second integrals at the rhs of \eqref{Lin1} and \eqref{Lin2}, we have to focus on the different intervals where $\alpha$ belongs.
\begin{itemize}
\item  \textbf{Case \ref{A1B}} By enlarging $s$, a double application of Young's inequality and bound \eqref{Cg} yield  
\begin{equation}\label{Lin3}
\const{lc} \int_\Omega |\nabla w|^{\frac{2(p+l)}{l}} \leq   \int_\Omega |\nabla w|^s + \const{le} \leq \const{lf} \quad \text{on } (0,\TM),
\end{equation}
and 
\begin{equation}\label{Lin4}
\const{ld} \int_\Omega |\nabla w|^{\frac{2(p-1)(p+l)}{p+l-2\alpha}} \leq   \int_\Omega |\nabla w|^s + \const{go} \leq \const{ho} \quad \text{for all } t \in (0,\TM).
\end{equation}
By plugging bounds \eqref{Lin1} and \eqref{Lin2} into \eqref{Est1L}, taking into account \eqref{Lin3} and \eqref{Lin4} and recalling that 
$\mu>2^{l} \chi \gamma \frac{(p-1)}{p}$, we can conclude reasoning as in Lemma \ref{Estim_general_For_u^p_nablaw^2rLemmaLog}.
\item  \textbf{Case \ref{A2B}} In this case, since $s$ cannot arbitrarily increase, the previous two integrals become, thanks to the Young inequality and $l>1$, for $\epsilon>0$ 
\begin{equation*}
\const{lc} \int_\Omega |\nabla w|^{\frac{2(p+l)}{l}} \leq \frac{p}{8(4p^2+n) \|w_0\|_{L^{\infty}(\Omega)}}  \int_\Omega |\nabla w|^{2(p+1)} + \const{li}  \quad \text{on } (0,\TM),
\end{equation*}
and 
\begin{equation}\label{Lin6}
\const{ld} \int_\Omega |\nabla w|^{\frac{2(p-1)(p+l)}{p+l-2\alpha}} \leq 
\frac{p}{8(4p^2+n) \|w_0\|_{L^{\infty}(\Omega)}} \int_\Omega |\nabla w|^{2(p+1)} + \const{ll} \quad \text{for all } t \in (0,\TM).
\end{equation}
The argument follows by invoking  \eqref{GradSfin} and $\mu>2^{l} \chi \gamma \frac{(p-1)}{p}$.
\item  \textbf{Case \ref{A3B}} First, we note that estimate \eqref{Lin6} still works taking $l=1$ and $\alpha<1$.
Moreover, with the aid of the Young inequality bound \eqref{Lin1} is turned into 
\begin{equation}\label{Lin1N}
\begin{split}
\frac{\xi^2 p(p-1)}{2} \int_\Omega (u+1)^p |\nabla w|^2 & \leq \frac{p}{8(4p^2+n) \|w_0\|_{L^{\infty}(\Omega)}^2} \int_\Omega |\nabla w|^{2(p+1)} \\ &
\quad + \xi^2 p^2 2^{\frac{2-p}{p}} \left(\frac{p-1}{p+1}\right)^{\frac{p+1}{p}} (4p^2+n)^{\frac{1}{p}}  \|\xi w_0\|_{L^{\infty}(\Omega)}^{\frac{2}{p}}
\int_\Omega (u+1)^{p+1} \quad \text{on } (0,\TM).
\end{split}
\end{equation}
Henceforth, by taking 
$$\mu > 2 \chi \gamma \frac{p-1}{p} +  \xi^2 p 2^{\frac{2}{p}} \left(\frac{p-1}{p+1}\right)^{\frac{p+1}{p}} (4p^2+n)^{\frac{1}{p}}  \|\xi w_0\|_{L^{\infty}(\Omega)}^{\frac{2}{p}},$$ we deduce the claim.
\item \textbf{Case \ref{A4B}} Finally, we consider $\alpha=l=1$. From the one hand, estimate \eqref{Lin1N} is still valid, while we have to deal in a different way with relation \eqref{Lin2}. By reasoning as in \cite[Lemma 5.2]{frassuviglialoro} and invoking \cite[(3.10) of Proposition 8]{YokotaEtAlNonCONVEX} with $s=1$, we find out that 
$\const{la}= p K^2 \|w_0\|_{L^{\infty}(\Omega)}^2 (n+p-1)$. Besides, an application of Young's inequality yields 
\begin{equation*}
\begin{split}
&p K^2 \|w_0\|_{L^{\infty}(\Omega)}^2 (n+p-1) \int_\Omega u^{2} |\nabla w|^{2p-2} \leq \frac{p}{8(4p^2+n) \|w_0\|_{L^{\infty}(\Omega)}^2}  \int_\Omega |\nabla w|^{2(p+1)}\\
&+p 2^{\frac{3p-1}{2}} K^{p+1} \left(\frac{n+p-1}{p+1}\right)^{\frac{p+1}{2}} ((p-1)(4p^2+n))^{\frac{p-1}{2}}  \|w_0\|_{L^{\infty}(\Omega)}^{2p}
\int_\Omega (u+1)^{p+1} \quad \text{for all } t \in (0,\TM).
\end{split}
\end{equation*}
Subsequently, for
\begin{equation*}
\begin{split}
\mu > 2 \chi \gamma \frac{p-1}{p} &+  \xi^2 p 2^{\frac{2}{p}} \left(\frac{p-1}{p+1}\right)^{\frac{p+1}{p}} (4p^2+n)^{\frac{1}{p}}  \|\xi w_0\|_{L^{\infty}(\Omega)}^{\frac{2}{p}} \\ & + 2^{\frac{3p+1}{2}} K^{p+1} \left(\frac{n+p-1}{p+1}\right)^{\frac{p+1}{2}} ((p-1)(4p^2+n))^{\frac{p-1}{2}}  \|w_0\|_{L^{\infty}(\Omega)}^{2p},
\end{split}
\end{equation*}
we can conclude.  
\end{itemize}
Continuity arguments in conjunction with the precise boundedness criterion mentioned in $\S$\ref{SectionLocalInTime} for the case $m_1=m_2=m_3=1$, allow us to conclude. 
\end{proof}
\begin{remark}\label{JustificationRemark}
Let us justify the second comment in Remark \ref{RemarkOnLinearCase}. The only term in \eqref{Estim1N} which has to be managed is (recall $n=2$, $l=m_1=m_2=m_3=1$) the integral
$(p-1)\chi \gamma \int_\Omega (u+1)^{p+1}$. By using the Gagliardo--Nirenberg inequality, we have
\begin{equation*}
\begin{split}
&(p-1)\chi \gamma \int_\Omega (u+1)^{p+1}= (p-1)\chi \gamma \|(u+1)^{\frac{p}{2}}\|_{L^{\frac{2(p+1)}{p}}(\Omega)}^{\frac{2(p+1)}{p}}\leq  (p-1)\chi\gamma C_{GN}\int_\Omega |\nabla (u+1)^{\frac{p}{2}}|^2+ \const{aoy} \quad \textrm{on } (0,\TM),
\end{split}
\end{equation*}
where $C_{GN}>0$ involves also the Gagliardo--Nirenberg constant. By inserting this bound into \eqref{Estim1N}, in order to proceed with the remaining steps, it is sufficient to make the term $\epsilon-(p-1)\left(\frac{4}{p}-\chi \gamma C_{GN}\right)<0$; this requires $\chi \gamma <\frac{4}{C_{GN}p}$ and continuity arguments, applied for $p>\frac{n}{2}=1$, make that it is also expressed as $\chi \gamma <\frac{4}{C_{GN}}.$ 
\end{remark}
\subsubsection*{\bf\textit{\quad Acknowledgments}}
The authors are very grateful to Professor Johannes Lankeit, who suggested to complement investigations \cite{FrassuLiViglialoro} and \cite{ChiyoFrassuViglialoro-Att-Rep-2022} with this article. 

SF and GV are members of the Gruppo Nazionale per l'Analisi Matematica, la Probabilit\`a e le loro Applicazioni (GNAMPA) of the Istituto Nazionale di Alta Matematica (INdAM), and are partially supported by the research projects \emph{Evolutive and Stationary Partial Differential Equations with a Focus on Biomathematics} (2019, Grant Number: F72F20000200007), {\em  Analysis of PDEs in connection with real phenomena} (2021, Grant Number: F73C22001130007), funded by  \href{https://www.fondazionedisardegna.it/}{Fondazione di Sardegna}. GV is also supported by MIUR (Italian Ministry of Education, University and Research) Prin 2017 \emph{Nonlinear Differential Problems via Variational, Topological and Set-valued Methods} (Grant Number: 2017AYM8XW). 
TL is partially supported by NNSF of P. R. China (Grant Number: 61503171), CPSF (Grant Number: 2015M582091), and NSF of Shandong Province (Grant Number: ZR2016JL021).

\end{document}